\documentclass[openany, amssymb, psamsfonts]{amsart}
\usepackage{mathrsfs,comment}
\usepackage[usenames,dvipsnames]{color}
\usepackage[normalem]{ulem}
\usepackage{url}
\usepackage[all,arc,2cell]{xy}
\UseAllTwocells
\usepackage{enumerate}
\usepackage{hyperref}
\usepackage{tikz}
\usetikzlibrary{calc}

\hypersetup{%
  bookmarksnumbered=true,%
  bookmarks=true,%
  colorlinks=true,%
  linkcolor=blue,%
  citecolor=blue,%
  filecolor=blue,%
  menucolor=blue,%
  pagecolor=blue,%
  urlcolor=blue,%
  pdfnewwindow=true,%
  pdfstartview=FitBH}

\def\makeautorefname#1#2{\expandafter\def\csname#1autorefname\endcsname{#2}}
\def\equationautorefname~#1\null{(#1)\null}
\makeautorefname{footnote}{footnote}%
\makeautorefname{item}{item}%
\makeautorefname{figure}{Figure}%
\makeautorefname{table}{Table}%
\makeautorefname{part}{Part}%
\makeautorefname{appendix}{Appendix}%
\makeautorefname{chapter}{Chapter}%
\makeautorefname{section}{Section}%
\makeautorefname{subsection}{Section}%
\makeautorefname{subsubsection}{Section}%
\makeautorefname{theorem}{Theorem}%
\makeautorefname{thm}{Theorem}%
\makeautorefname{cor}{Corollary}%
\makeautorefname{lem}{Lemma}%
\makeautorefname{prop}{Proposition}%
\makeautorefname{pro}{Property}
\makeautorefname{conj}{Conjecture}%
\makeautorefname{defn}{Definition}%
\makeautorefname{notn}{Notation}
\makeautorefname{notns}{Notations}
\makeautorefname{rem}{Remark}%
\makeautorefname{quest}{Question}%
\makeautorefname{exmp}{Example}%
\makeautorefname{ax}{Axiom}%
\makeautorefname{claim}{Claim}%
\makeautorefname{ass}{Assumption}%
\makeautorefname{asss}{Assumptions}%
\makeautorefname{con}{Construction}%
\makeautorefname{prob}{Problem}%
\makeautorefname{warn}{Warning}%
\makeautorefname{obs}{Observation}%
\makeautorefname{conv}{Convention}%

\newtheorem{cor}{Corollary}[section]
\newtheorem{prop}{Proposition}[section]
\newtheorem{lem}{Lemma}[section]

\theoremstyle{definition}
\newtheorem{defn}{Definition}[section]

\newtheorem{exmp}{Example}[section]

\newtheorem{notns}{Notations}[section]

\newtheorem{rem}{Remark}[section]

\makeatletter
\let\c@obs=\c@thm
\let\c@cor=\c@thm
\let\c@prop=\c@thm
\let\c@lem=\c@thm
\let\c@prob=\c@thm
\let\c@con=\c@thm
\let\c@conj=\c@thm
\let\c@defn=\c@thm
\let\c@notn=\c@thm
\let\c@notns=\c@thm
\let\c@exmp=\c@thm
\let\c@ax=\c@thm
\let\c@pro=\c@thm
\let\c@ass=\c@thm
\let\c@warn=\c@thm
\let\c@rem=\c@thm
\let\c@sch=\c@thm
\let\c@equation\c@thm
\numberwithin{equation}{section}
\makeatother

\title{A Study in Markov Chains, Loop-Erased Random Walk and Loop Soups}

\author{Zhuohan (Joshua) Gu}

\date{September, 2023}

\begin{document}

\begin{abstract}

In this paper, we make a few random explorations that relate directly to the items mentioned in the title. We define transient chains and recurrent chains with ``killing'', the Green's function, the Laplacian operator, and harmonic functions. We then introduce the loop-erased random walk (LERW) and its relationship with the uniform spanning tree (UST). We finish by introducing loop measures and soups and defining the ``growing  loop'' model.

\end{abstract}

\maketitle

\tableofcontents

\section{Introduction}

This paper takes the first three chapters of the book by Gregory F. Lawler \cite{ams_RE} and summarizes them. Instead of presenting the proofs in the book, we present original proofs for important lemmas, propositions, and facts in the book that are given without proofs, as well as solutions to interesting challenging problems in the book that are given without solutions. Missing definitions and missing proofs of theorems, lemmas, propositions, corollaries, etc., can be found in the book \cite{ams_RE}.

Chapter 2 mainly focuses on discrete time chains, but we also discuss how to construct a continuous time Markov chain later in the chapter. Analysis of the loop-erasing procedure and its applications are introduced in Chapter 3. Chapter 4 introduces loop measures, loop soups and the ``growing looop'' model. Three different ways of constructing the loop configurations are introduced: (ordered) growing loops, rooted loops, and unrooted loops. Chapter 3 (Loop-Erased Random Walk) and Chapter 4 (Loop Soups) depend on the material in Chapter 2 (Markov Chains). Chapter 4 uses Chapter 3.

\section{Markov Chains}

\subsection{Definition}
\hfill\\

If $A$ is finite with $N$ elements, the \textit{transition probabilities} $p(x, y)$ of a \textit{finite Markov chain} on $A$ are often collected into a $N \times N$ matrix $P = [p(x, y)]$ called the \textit{transition matrix}.

The \textit{n-step transition matrix} is obtained by raising the \textit{transition matrix} $P$ to the \textit{n}th power. One important property of a transition matrix is that each of the rows sum to 1, that is,
\begin{equation}\label{(2.4)}
    \sum_{y \in A} p(x, y) = 1.
\end{equation}

\begin{prop}
    (Chapman-Kolmogorov equations) If $n, m \in \mathbb{N}$ and $x, y \in A,$
    \[ p_{n+m}(x, y) = \sum_{z \in A}p_n(x, z)p_m(z, y). \]
\end{prop}
\begin{proof}
    The proof uses the law of total probability and the definition of conditional probability. For a complete proof, see the proof of Proposition 1.1 in \cite{ams_RE}.
\end{proof}
The \textit{Chapman-Kolmogorov equations} can be written in matrix form as
\[ P^{n+m} = P^n \cdot P^m, \]
where $\cdot$ denotes matrix multiplication.

\begin{defn}
    A Markov chain is \textit{irreducible} if for every $x, y \in A$ there exists an integer $n\geq0$ such that $p_n(x, y) > 0$.
\end{defn}
That is, it is possible to get from any state to any other state. 

\begin{defn}
     Let $x \in A$ and $\tau_x = \min \{ k \geq 1: X_k = x\}$. $\tau_x$ is the index of the first visit to $x$ after time $0$. If $\mathbb{P}\{\tau_x < \infty \mid X_0 = y\} = 1$, then the irreducible Markov chain is \textit{recurrent}. If $\mathbb{P}\{\tau_x < \infty \mid X_0 = y\} < 1$, then the irreducible Markov chain is \textit{transient}.
\end{defn}
In other words, an irreducible Markov chain is \textit{recurrent} if eventual return to every point is certain.

\begin{notns}
    If E is an event and Y is a random variable, then
    \[\mathbb{P}^x(E) = \mathbb{P}\{E \mid X_0 = x\},\]
    \[\mathbb{E}^x(Y) = \mathbb{E}\{Y \mid X_0 = x\}.\]
\end{notns}

\begin{notns} (Indicator function notation)
    If $E$ is an event, then $1_E$ is the random variable that equals 1 if $E$ occurs and equals 0 if $E$ does not occur.
\end{notns}

We can use the indicator function notation $1_E$ to denote the total number of visits to the state $x$,
\[V_x = \sum_{n=0}^{\infty} 1\{X_n = x\}.\]

\begin{lem} If $V_y$ denotes the total number of visits to the state $y$, then
    \[\mathbb{E}^x[V_y] = \sum_{n=0}^{\infty} p_n(x, y).\]
\end{lem}
\begin{proof} Note that $\mathbb{E}[1_E] = \sum_{x \in E} 1 \cdot \mathbb{P}\{X=x\} = \mathbb{P}(E)$ and hence
    \begin{equation*}
        \mathbb{E}[V_y \mid X_0 = x] = \sum_{n=0}^{\infty}\mathbb{E}[1\{X_n = y\} \mid X_0 = x] = \sum_{n=0}^{\infty} p_n(x, y).
    \end{equation*}
\end{proof}

\begin{rem}
    In this paper, the word ``visit'' includes the visit at time 0 but the word ``return'' refers only to the visits after time 0. Therefore, $\tau_x$ describes the time it takes for a Markov chain to return to state $x$, and $V_x$ is the total number of visits to $x$. 
\end{rem}

\begin{prop} \label{prop2.10}
    If $X_n$ is an irreducible Markov chain and $x \in A$, then
    \[\mathbb{P}\{\tau_x = \infty\} = \frac{1}{\mathbb{E}^x[V_x]}.\]
    In particular, $\mathbb{P}\{\tau_x = \infty\} = 0$ if and only if $\mathbb{E}^x[V_x] = \infty$.
\end{prop}

\begin{prop} \label{prop2.11}
    The following are equivalent for an irreducible Markov chain.
\begin{enumerate}
    \item The chain is recurrent.
    \item For every $x, y \in A$, $\mathbb{E}^x[V_y] = \infty$.
    \item There exists $x, y \in A$ such that $\mathbb{E}^x[V_y] = \infty$.
    \item For every $x, y \in A$, $\mathbb{P}^x[V_y = \infty] = 1$.
    \item There exists $x, y \in A$ such that $\mathbb{P}^x[V_y = \infty] = 1$. 
\end{enumerate}
\end{prop}

\begin{proof}
    We will prove the equivalence by proving: $(1) \Longrightarrow (2) \Longrightarrow (3) \Longrightarrow (4) \Longrightarrow (5) \Longrightarrow (1)$.
    
    $(1) \Longrightarrow (2):$ Suppose the chain is recurrent. By definition, for all $x, y \in A$, $\mathbb{P}\{\tau_y < \infty \mid X_0 = x\} = 1$. In other words, starting in state $x$, the chain will certainly return to $y$ in finite steps. Since the chain is recurrent, $\mathbb{P}^{y}\{\tau_y = \infty \} = 0$. By \autoref{prop2.10}, this implies $\mathbb{E}^y[V_y] = \infty$ and hence $\mathbb{E}^x[V_y] = \infty$.
    
    $(2) \Longrightarrow (3):$ Clear, since $(3)$ is a case of $(2)$.

    $(3) \Longrightarrow (4):$ Suppose there exists $x, y \in A$ such that $\mathbb{E}^x[V_y] = \infty$. Then, $\mathbb{E}^y[V_y] = \infty$. This implies that starting at $y$, the chain will certainly return to $y$ in finite steps. If not, it contradicts $\mathbb{E}^y[V_y] = \infty$. We can perceive this as having infinitely many finite-step ``loops'' rooted at $y$, where each ``loop'' only visits $y$ at the root. For all $z, w \in A$, since there exists $n\geq0$ such that $p_n(y, z)  > 0$ by irreducibility, $z$ must be in one of the ``loops'' and hence $\mathbb{P}^z\{\tau_y < \infty \} = 1$. Since there are infinitely many such loops and $p_n(y, w) > 0$, $\mathbb{P}^z\{V_w = \infty \} = 1$.
    
    $(4) \Longrightarrow (5):$ Clear, since $(5)$ is a case of $(4)$.

    $(5) \Longrightarrow (1):$ Suppose there exists $x, y \in A$ such that $\mathbb{P}^x\{V_y = \infty\} = 1$. Then, $\mathbb{P}^y\{V_y = \infty\} = 1$ so there are infinitely many such ``loops'' as defined previously. For all $z, w \in A$, since there exists $n\geq0$ such that $p_n(y, z)  > 0$, $z$ must be in one of the ``loops'' and hence $\mathbb{P}\{\tau_y < \infty \mid X_0 = z\} = 1$. Since there are infinitely many such loops and $p_n(y, w) > 0$, the chain will reach $w$ from $y$ in finite steps and therefore $\mathbb{P}\{\tau_w < \infty \mid X_0 = z\} = 1$.
\end{proof}

\begin{exmp} \label{Ex1.4}
    Consider the Markov chain whose state space is the integers with transition probabilities
    \[p(x, x+1) = q, p(x, x-1) = 1-q,\]
    where $0 < q < 1$. Let
    \[p_n = p_n(0, 0) = \mathbb{P}^0\{X_n = 0\}.\]

    (1) $P_n = 0$ if $n$ is odd since we need an even number of steps to get back to the starting point.

    (2) $p_{2n} = \binom{2n}{n}q^n(1-q)^n$ since we need $n$ steps of ``+1'' and $n$ steps of ``-1'' to get back to the starting point.

    (3) If $q = \frac{1}{2}$, assuming Stirling's formula, we get
    \[\lim_{n \rightarrow \infty} n^{1/2}p_{2n} = \frac{1}{\sqrt{\pi}}.\]
    (4) Suppose $q=1/2$. By (3), for all $\epsilon > 0$, there exists $N > 0$ such that for all $n>N$ we have
    \[\frac{1}{\sqrt{\pi}}-\epsilon < n^{1/2}p_{2n} < \frac{1}{\sqrt{\pi}} + \epsilon.\]
    Then, \[\mathbb{E}^0[V_0] = \sum_{n=0}^{\infty}p_{2n}(0,0)>\sum_{n>N}^{\infty}(\frac{1}{\sqrt{\pi}}-\epsilon)\frac{1}{\sqrt{n}} \rightarrow \infty.\]
    By \autoref{prop2.11}, the chain is recurrent.

    (5)Suppose $q \neq 1/2$:
    \[\frac{p_{2n, q \neq 1/2}}{p_{2n, q = 1/2}} = \frac{\binom{2n}{n}q^n(1-q)^n}{\binom{2n}{n}(1/2)^{2n}} = [4q(1-q)]^n.\]
    Then we get
    \[p_{2n, q \neq 1/2} = p_{2n, q=1/2}\cdot[4q(1-q)]^n.\]
    Similarly,
    \[\mathbb{E}^0[V_0] = \sum_{n=0}^{\infty}p_{2n} \leq \sum_{n=0}^{N}p_{2n} + (\frac{1}{\sqrt{\pi}}+\epsilon)\sum_{n>N}^{\infty}\frac{1}{\sqrt{n}}[4q(1-q)]^n,\]
    which converges using the ratio test. Therefore, the chain is transient for $q \neq 1/2$.
\end{exmp}

\subsection{Laplacian and harmonic functions}
\begin{defn}
    The \textit{Laplacian} is the linear transformation $\mathcal{L} = I - P$ where $I$ denotes the identity matrix and $P$ is the transition matrix,
    \[\mathcal{L}f(x) = (I-P)f(x) = f(x) - \sum_{y \in A} p(x, y)f(y).\]
\end{defn}
Recalling \autoref{(2.4)}, $\sum_{y \in A} p(x, y)f(y)$ can be seen as the ``mean value" or ``expected value" of $f$ along $x$'s adjacent vertices.

For a random walk on a simple graph, we can write
\[\mathcal{L}f(x) = f(x) - \mathbb{E}^x[f(X_1)] = f(x) - MV(f; x)\] where $\mathbb{E}^x[f(X_1)]$ denotes the expected value of $f$ at time $1$ given $X_0 = x.$ $X_1$ is a random variable taking values from $x$'s neighbors. $MV(f;x)$ denotes the ``mean value'' of $f$ along $x$'s neighbors.

If $f : \mathbb{R}^d \rightarrow \mathbb{R}$ is a smooth function, the Laplacian is defined as
\[\Delta f(x) = \sum_{j=1}^{d} \partial^2_{x_jx_j}f(x).\]

\begin{prop}
    Suppose $D \subset \mathbb{R}^d$ is open, and $f : D \rightarrow \mathbb{R}^d$ is a function with continuous first and second derivatives. Then, for $x \in D$,
    \[\Delta f(x) = 2d \lim_{\epsilon \downarrow 0} \frac{MV(f; x, \epsilon)-f(x)}{\epsilon^2},\]
    where $MV(f; x, \epsilon)$ denotes the mean (average) value of $f$ on the sphere of radius $\epsilon$ centered at $x$, $\{y \in \mathbb{R}^d : |x-y| = \epsilon\}$.
\end{prop}

\begin{proof}
    It is the equivalent of showing \[\Delta f(\vec{0}) = 2d \lim_{\epsilon \downarrow \vec{0}} \frac{MV(f; \vec{0}, \epsilon)-f(\vec{0})}{\epsilon^2}.\]
    Expand $f$ around $\vec{0}$ in its Taylor polynomial of degree 2:
    \[f(\vec{x}) = f(\vec{0}) + \sum_{i=1}^{d}\frac{\partial f}{\partial x_i}(\vec{0})x_i + \frac{1}{2}\sum_{i,j=1}^{d}\frac{\partial^2 f}{\partial x_i x_j}(\vec{0})x_ix_j + o(|x|^2).\]

    Substitute the expansion into the right-hand side:
    \begin{equation*}
        \begin{split}
            & 2d\lim_{\epsilon \downarrow \vec{0}}\frac{1}{\epsilon^2}\left[\frac{1}{z_\epsilon}\int_{|x|=\epsilon} f(\vec{0}) + \sum_{i=1}^{d}\frac{\partial f}{\partial x_i}(\vec{0})x_i + \frac{1}{2}\sum_{i,j=1}^{d}\frac{\partial^2 f}{\partial x_i x_j}(\vec{0})x_ix_j + o(|x|^2)\,ds(\vec{x}) - f(\vec{0})\right]\\
            & = d\lim_{\epsilon \downarrow \vec{0}}\frac{1}{\epsilon^2}\left[\frac{1}{z_\epsilon}\int_{|x|=\epsilon} \sum_{i=1}^{d}\frac{\partial^2 f}{\partial x_i^2}(\vec{0})x_i^2 \,ds(\vec{x})\right]\\
            & = d\lim_{\epsilon \downarrow \vec{0}}\frac{1}{\epsilon^2}\left[\sum_{i=1}^{d}\frac{\partial^2 f}{\partial x_i^2}(\vec{0})\frac{\epsilon^2}{d}\right]\\
            & = \Delta f(\vec{0}).
        \end{split}
    \end{equation*}
    $z_\epsilon$ is the ``surface area'' of the sphere. The first equality uses the symmetry in the integral. Note that
    \[\int_{|x|=\epsilon} x_1^2 \,ds(\vec{x}) = \int_{|x|=\epsilon} x_2^2 \,ds(\vec{x}) = \cdots = \int_{|x|=\epsilon} x_d^2 \,ds(\vec{x}).\]
    Then we can see that \[\int_{|x|=\epsilon} x_i^2 \,ds(\vec{x}) = \frac{1}{d} \int_{|x|=\epsilon} x_1^2 + x_2^2 + \cdots + x_d^2 \,ds(\vec{x}) = \frac{\epsilon^2}{d},\]
    which is used in the second equality.
\end{proof}

\begin{defn}
    A function $f$ is \textit{(P-)harmonic} on $A' \subset A$ if for all $x \in A'$, $Pf(x) = f(x)$, that is
    \[\mathcal{L}f(x) = 0, \text{  } x \in A'.\] 
\end{defn}

\subsection{Markov chain with boundary}
\hfill\\

A state space $\overline{A}$ can be written as
\[\overline{A} = A \cup \partial A\]
where $A$ denotes the \textit{interior vertices} and $\partial A$ are the \textit{boundary vertices}. Here, the usages of $A, \partial A, \overline{A}$ are directly analogous to a set, its limit points, and its closure, respectively.

If $P$ is a transition matrix for an irreducible Markov chain on $\overline{A}$ with entries $p_{ij}$ for $x_i, x_j \in \overline{A}$, then $\tilde{P}$ is the corresponding matrix such that for each $x_i \in \partial A$, we change the entries in the $i^{th}$-row:
\[\tilde{p}(x_i, x_j) = \begin{cases}
    1 & \text{if $x_i = x_j$}\\
    0 & \text{if $x_i \neq x_j$}
  \end{cases}.\]
We let $P_A$ be the submatrix of $P$ obtained by restricting the states in A, that is, for all $P_{ij}$ in $P_A$, $x_i$ and $x_j$ are in $A$. This is the same as $\tilde{P}_A$ and $P_A = \tilde{P}_A$.

\begin{defn}
    Given a stochastic process $\{X_0, X_1, X_2, ...\}$, a non-negative integer random variable $\tau$ is called a \textit{stopping time} if for all integers $k \geq 0$, $\tau \leq k$ depends only on $X_0, X_1, ..., X_k$.
\end{defn}
Assume $\partial A \neq \emptyset$ and define the stopping time
\[T = T_A = \min\{k \geq 0: X_k \notin A\}.\]
In other words, $T_A$ is the first index to leave $A$, and $T_A$ is a random variable. If $A$ is finite, since the chain is irreducible,
\[\mathbb{P}\{T < \infty\} = 1.\]

\begin{defn}
    The \textit{Poisson kernel} is the function $H_A: \overline{A} \times \partial A \rightarrow [0, 1]$ given by
    \[H_A(x, z) = \mathbb{P}^x\{X_T = z\} = \mathbb{P}\{X_T = z \mid X_0 = x\}.\]
\end{defn}
The Poisson kernel describes the probability of leaving $A$ at a particular point in $\partial A$ starting at $x$. Note that for each $x$, 
\[\sum_{z \in \partial A} H_A(x, z) = 1.\]

If $x \in \partial A$, that is, the chain that starts outside of $A$, then $T = 0$ and hence
\begin{equation}\label{(2.21)}
    H_A(x, z) = \begin{cases}
        1, & \text{$x = z$}\\
        0, & \text{$x \neq z$}
        \end{cases},\text{ } x \in \partial A.
\end{equation}

\begin{prop}\label{prop2.22}
    Suppose $P$ is an irreducible transition matrix on $\overline{A} = A \cup \partial A$. If $z \in \partial A$ and $h(x) = H_A(x, z)$, then h is the unique bounded function on $\overline{A}$ that is harmonic in A and satisfies the boundary condition \autoref{(2.21)} on $\partial A$.
\end{prop}

 \autoref{prop2.22} establishes the existence and uniqueness of bounded harmonic functions on Markov chains satisfying specific \textit{Dirichlet boundary condition}. If the proposition generalizes to arbitrary boundary conditions, it is called the solution to the \textit{Dirichlet problem}.

\begin{prop}
    (Bounded Convergence Theorem) Suppose $X_1, X_2, ...$ is a collection of random variables such that with probability one, the limit
    \[X = \lim_{n \rightarrow \infty} X_n\]
    exists. Assume also that there exists $J < \infty$ such that $|X_n| \leq J$ for all $n$. Then
    \[\mathbb{E}[X] = \lim_{n \rightarrow \infty} \mathbb{E}[X_n].\]
\end{prop}

\begin{proof}
    We will first show that it suffices to prove this when $X$ is identically equal to 0. Then, for n sufficiently large, $\mathbb{E}[|X_n|] \leq 2\epsilon$.

    Let $Y_n = X_n - X$. Then
    \[\lim_{n \rightarrow \infty} Y_n = \lim_{n \rightarrow \infty} (X_n - X).\]
    Since $X_n \leq J$ for all $n$,  $Y_n \leq 2J$ for all $n$. Then
    \[\lim_{n \rightarrow \infty} \mathbb{E}[Y_n] = \lim_{n \rightarrow \infty} \mathbb{E}[X_n - X] =  \lim_{n \rightarrow \infty} \mathbb{E}[X_n] - \mathbb{E}[X] = 0\]
    and hence
    \[\mathbb{E}[X] = \lim_{n \rightarrow \infty} \mathbb{E}[X_n].\]
    To show $\mathbb{E}[|X_n|] \leq 2\epsilon$, for every $\epsilon > 0$, we write,
    \begin{equation*}
        \begin{split}
            \mathbb{E}[|X_n|] & = \mathbb{E}[|X_n|\cdot 1\{|X_n| \leq \epsilon\}] + \mathbb{E}[|X_n|\cdot 1\{|X_n| > \epsilon\}]\\
            & \leq \epsilon + \mathbb{E}[J \cdot 1\{|X_n| > \epsilon\}]\\
            & = \epsilon + J \cdot \mathbb{P}\{|X_n| > \epsilon\}.
        \end{split}
    \end{equation*}
    Since $\mathbb{P}\{|X_n| > \epsilon\} \rightarrow 0$ as $n \rightarrow \infty$, there exists $N$ such that for all $n > N$, $\mathbb{P}\{|X_n| > \epsilon\} \leq \frac{\epsilon}{J}$. Then
    \[\mathbb{E}[|X_n|] \leq \epsilon + J \cdot \frac{\epsilon}{J} = 2\epsilon.\]
\end{proof}

\subsection{Green's function}
\begin{defn}
    Assume either $\partial A \neq \emptyset$ or the chain is transient. The \textit{Green's function} $G_A(x, y)$ is defined for $x, y \in A$ by
    \[G_A(x, y) = \mathbb{E}^x[V_y] = \sum_{n=0}^{\infty}p_n(x, y).\]
\end{defn}
The Green's function describes the expected number of visits to $y \in A$ starting at $x \in A$. We also consider $G_A$ as a linear transformation
\[G_Af(x) = \sum_{y \in A}G_A(x, y)f(y) = \sum_{y\in A}\sum_{n=0}^{\infty}p_n(x, y)f(y).\] In other words,
\begin{equation}
    G_A = I + P_A + P^2_A + P^3_A + \cdots = (I-P_A)^{-1} = \mathcal{L}^{-1}_A.
\end{equation}

\begin{exmp}
    Consider a binary tree defined as follows. Let $A$ be the set of finite sequences of 0s and 1s such as 0010110. We include the empty sequence which we represent as $\emptyset$. We say that sequence $a$ is the parent of sequence $b$ if $b$ is of the form $ar$ where $r$ is 0 or 1. All sites have a parent except for the empty sequence. Consider the Markov chain with transition probabilities
    \[p(a, b) = \frac{1}3{}, \quad p(b, a) = \frac{1}{3},\]
    if $a$ is a parent of $b$. In other words, the random walker chooses randomly among its parent and its two children. Since the empty sequence $\emptyset$ has no parent we also set $p(\emptyset, \emptyset) = \frac{1}{3}$.

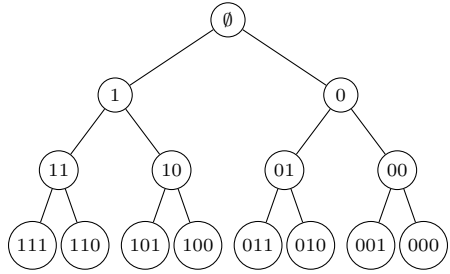
\begin{center}
    \begin{tikzpicture}[
        level 1/.style={level distance=10mm,sibling distance=30mm},
        level 2/.style={level distance=10mm,sibling distance=15mm},
        level 3/.style={level distance=10mm,sibling distance=7mm},
        font=\scriptsize,inner sep=2pt,every node/.style={draw,circle,minimum size=3ex}]
        
        \node {$\emptyset$} 
        child {node {1} 
                child {node{11} 
                    child {node{111}} 
                    child {node{110}}
                } 
                child {node{10} 
                    child {node{101}}
                    child {node{100}}
                }
            edge from parent }
        child {node {0} 
            child {node {01} 
                child {node{011}}
                child {node{010}}
            }
            child{node {00} 
                child {node{001}}
                child {node{000}}
            }
                } ;
    \end{tikzpicture} 
\end{center}

(1) Observe that the probability of going up is $1/3$ and the probability of going down is $2/3$. By \autoref{Ex1.4}, $q \neq 1/2$ and hence the chain is transient.

(2)For $b \neq \emptyset$, let $p(b)$ be the probability that the chain starting at $b$ ever reaches the empty sequence and let $|b|$ be the length of (number of digits in) $b$. Suppose the probability that $b$ ever reaches its parent is $\lambda$. Then the probability that $b$'s parent ever reaches its parent is also $\lambda$. Then $p(b)$ must be of the form $\lambda^{|b|}$.

(3) From (2) it follows that
\[\lambda = \frac{1}{3} + \frac{2}{3}\lambda^2.\]
By solving this equation we find that $\lambda = 1 \text{ or } 1/2$. Since $\lambda \neq 1$, $\lambda = 1/2$.

(4) To find $G(\emptyset, \emptyset)$, recall that $G(\emptyset, \emptyset) = \mathbb{E}^\emptyset[V_\emptyset].$ Note that $V_\emptyset$ denotes the number of visits to $\emptyset$. In other words, $V_\emptyset$ denotes the number of visits to $\emptyset$ until the chain never returns to $\emptyset$. Then $V_\emptyset$ has a geometric distribution, written $V_\emptyset \sim Geo(1/2)$, and hence
\[G(\emptyset, \emptyset) = \frac{1}{1/2} = 2.\]

(5) It is also interesting to find $G(b, \emptyset)$ and $G(\emptyset, b)$ if $b$ is any sequence. Note that $G(b, \emptyset) = p(b)G(\emptyset, \emptyset)$. Then \[G(b, \emptyset) = \lambda^{|b|} \cdot 2 = 2 \cdot (\frac{1}{2})^{|b|}.\]
To find $G(\emptyset, b)$, suppose $|b| = 2$ and call it level 2. Let $L2$ denote level 2. Observe that
\[G(\emptyset, 11) = G(\emptyset, 10) = G(\emptyset, 01) = G(\emptyset, 00) = \frac{1}{4}G(\emptyset, L2) = \frac{1}{4}G(L2, L2),\]
where $G(\emptyset, L2)$ denotes the expected visits to level 2 starting at $\emptyset$ and $G(L2, L2)$ denotes the expected visits to level 2 starting at level 2. The last equality holds because starting at $\emptyset$, the chain will certainly reach $L2$ by transient. Recall in (4) that $V_{L2}$ has a geometric distribution and
\[\mathbb{P}^{L2}\{\text{never returns to $L$2}\} = \frac{2}{3} \cdot \frac{1}{2} = 1/3.\] Then $G(L2, L2) = \frac{1}{1/3} = 3$ and hence
\[G(\emptyset, b) = \frac{3}{2^{|b|}}.\]

(6) Suppose $b \neq \emptyset$. Define $\tilde{p}(b)$ to be the probability that $b$ ever reaches the other side of the binary tree. Then $\tilde{p}(b)$ is a bounded nonconstant function on $A$ that is harmonic with respect to this chain.
\end{exmp}

\subsection{Continuous time}
\hfill\\

We have been discussing Markov chains indexed by integer times. We now consider continuous-time process, Markov chains $Y_t$ indexed by time $t \in [0, \infty)$. These can be constructed from discrete-time chains by assuming that each time the process reaches a state, the amount of time it spent at the site before taking the next step has an \textit{exponential distribution}.

The exponential distribution satisfies the memoryless property:
\[\mathbb{P}\{T>t+s \mid T>s\} = \mathbb{P}\{T>t\},\]
which is important in order to construct a Markov chain. It tells us that the probability of leaving a site soon does not depend on how long one has been at the site.

\begin{prop}
    Suppose $T$ is a nonnegative random variable with a continuous distribution function $F$ satisfying the memoryless property: for all $t, s > 0$,
    \[\mathbb{P}\{T>t+s \mid T>s\} = \mathbb{P}\{T>t\}.\]
    Then $T$ has an exponential distribution.
\end{prop}

\begin{proof}
    Define $f(x) = \log \mathbb{P}\{T>x\}.$ For all $t, s > 0$,
    \begin{equation*}
        \begin{split}
            f(t+s) & = \log \mathbb{P}\{T>t+s\}\\
            & = \log \mathbb{P}\{T > s \cap T > t+s\}\\
            & = \log [\mathbb{P}\{T>s\}\mathbb{P}\{T>t+s \mid T>s\}]\\
            & = \log [\mathbb{P}\{T>s\}\mathbb{P}\{T>t\}]\\
            & = f(s) + f(t).
        \end{split}
    \end{equation*}
    By the continuity of $f$, $f(x) = xf(1)$ for all $x > 0$. Since $T$ is nonnegative, $f(0) = 1$. Note that $f(1) < 0$. Set $f(1) = -\lambda$ for some $\lambda > 0$. Then $f(x) = -\lambda x = \log \mathbb{P}\{T>x\}$. Then we get $\mathbb{P}\{T>x\} = e^{-\lambda x}$ and hence
    \[F(x) = \mathbb{P}\{T \leq x\} = 1 - \mathbb{P}\{T>x\} = 1 - e^{-\lambda x},\]
    which is the cumulative distribution function of an exponential distribution.
\end{proof}

\begin{prop}
    Suppose $T_1, T_2, ..., T_n$ are independent exponential random variables with rates $\lambda_1, \lambda_2, ..., \lambda_n$. Let $T = \min\{T_1, T_2, ..., T_n\}$. Then $T$ has an exponential distribution with rate \[\sum_{i=1}^{n} \lambda_i.\]
\end{prop}

\begin{proof}
    We will find the cumulative distribution function of $T$, written as $F(t) = \mathbb{P}\{T \leq t\}$:
    \begin{equation*}
        \begin{split}
            F(t) & = 1 - \mathbb{P}\{T_i > t \text{ for all } i\}\\
            & = 1 - \prod_{i=1}^{n} \mathbb{P}\{T_i > t\}\\
            & = 1 - \prod_{i=1}^{n} e^{- \lambda_i t}\\
            & = 1 - e^{- \sum_{i=1}^{n}\lambda_i t}.
        \end{split}
    \end{equation*}
\end{proof}

The exponential distribution of $T$ is used in the construction of the continuous time chain. Whenever we reach a state $x$, we wait an exponential amount of time with rate $\lambda(x)$, and then we move to state $y$ with rate $\lambda(x, y) = \lambda(x)p(x, y)$. $T_i$ can be seen to be the waiting time until the chain moves to $y_i$, written as $T_i \sim exp(\lambda(x, y_i))$. Then $T$ denotes the first ``neighbor'' of $x$ to call $x$ and say, ``jump here!''

\section{Loop-Erased Random Walk}

In this chapter, we define and explore the loop-erased random walk (LERW) and show how the LERW relates to the uniform spanning trees (UST).

\subsection{Loop erasure}
\begin{defn}
    Write a path $\omega$ of length $n$ as a finite sequence of points
    \begin{equation}
        \omega = [\omega_0, \ldots, \omega_n]
    \end{equation}
    with $\omega_j \in \overline{A}$. We call a path a \textit{self-avoiding walk} (SAW) of length $n$ if all of the vertices $\{\omega_0, \ldots, \omega_n\}$ are distinct.
\end{defn}

\begin{defn}
    If $\omega = [\omega_0, \omega_1, \ldots, \omega_n]$ is a path, then its \textit{(chronological) loop erasure}, denoted by $LE(\omega)$, is the SAW $\eta = [\eta_0, \eta_1, \ldots, \eta_k]$ defined as follows.
    \begin{itemize}
        \item Set $\eta_0: = \omega_0$.
        \item If $\omega_n = \omega_0$, we set $k=0$ and terminate; otherwise, let $\eta_1$ be the first vertex in $\omega$ after the last visit to $\omega_0$, that is, $\eta_1: = \omega_{i+1}$, where $i: = max\{j; \omega_j = \omega_0\}$.
        \item If $\omega_n = \eta_1$, then we set $k=1$ and terminate; otherwise, let $\eta_2$ be the first vertex in $\omega$ after the last visit to $\eta_1$, and so on.
    \end{itemize}
\end{defn}

$LE(\omega)$ is obtained by erasing cycles in $\omega$ in the order they appear. Note that the definition of $LE(\omega)$ depends on the order in which the vertices $[\omega_0, \omega_1, \ldots, \omega_n]$ are traversed. This procedure is sometimes called ``forward loop-erasing'' to distinguish it from ``backward loop-erasing'', which is defined as follows.
\begin{defn}
    Suppose $\omega = [\omega_0, \omega_1, \ldots, \omega_n]$. Let $\omega^R = [\omega_n \omega_{n-1}, \ldots, \omega_0]$ be the path $\omega$ traversed in the opposite direction. Then
    \[\tilde{\eta} = LE(\omega^R), \quad LE^R(\omega) = \tilde{\eta}^R.\]
\end{defn}

\begin{exmp}
    The following is an example of a path $\omega$ for which $LE^R(\omega) \neq LE(\omega)$.

    Suppose $\omega = [\omega_0, \omega_1, ..., \omega_7]$, which $\omega_1 = \omega_4$ and $\omega_2 = \omega_6$. Then $LE(\omega) = [\omega_0, \omega_4, \omega_5, \omega_6, \omega_7]$ and $LE^R(\omega) = [\omega_0, \omega_1, \omega_2, \omega_7]$. Hence $LE^R(\omega) \neq LE(\omega)$.
\end{exmp}

\subsection{Loop-erased random walk}
\hfill\\

In this section, we investigate the ``loop-erased random walk (LERW) from $x \in A$ to $\partial A$'' where by this we mean start the Markov chain in state $x$, stop it when it first reaches $\partial A$.

If $x \in A$, let $\mathcal{K}_A(x, \partial A)$ denote the set of paths starting at $x$ stopped at the first time that they reach the boundary. We write
\[\mathcal{K}_A(x, \partial A) = \bigcup_{z \in \partial A} \mathcal{K}_A(x, z)\]
where $\mathcal{K}_A(x, z)$ denotes the set of such paths that end at $z \in \partial A$. The probability measure $p$ on $\mathcal{K}_A(x, \partial A)$ is $p[\mathcal{K}_A(x, z)] = H_A(x, z)$: if $\omega = [\omega_0, \ldots, \omega_n] \in \mathcal{K}_A(x, \partial A),$
\[p(\omega) = \mathbb{P}^x\{X_1=\omega_1, X_2=\omega_2, \ldots, X_n=\omega_n\}=\prod_{j=1}^{n}p(\omega_{j-1}, \omega_j).\]

Let $\mathcal{R}_A(x, \partial A)$ denote the set of self-avoiding walks (SAWs) in $\mathcal{K}_A(x, \partial A)$ and define $\mathcal{R}_A(x, z)$ similarly.

\begin{defn}
    \textit{Loop-erased random walk (LERW)} from $x$ to $\partial A$ is the probability measure \[\hat{p}(\eta) = \hat{p}_A(\eta) = \sum_{\omega \in \mathcal{K}_A(x, \partial A), LE(\omega)=\eta}p(\omega).\]
\end{defn}

Since the initial and terminal vertices are fixed in the loop-erasing
procedure we see that for all $z \in \partial A$,
\[\hat{p}[\mathcal{R}_A(x, z)] = p[\mathcal{K}_A(x, z)] = H_A(x, z).\]

\begin{prop} \label{prop3.7}
    Suppose that $\eta = [\eta_0, \eta_1, \ldots, \eta_k] \in \mathcal{R}_A(x, \partial A)$ and let $A_j = A \backslash \{\eta_0, \ldots, \eta_{j-1}\}$. Then
    \[\hat{p}_A(\eta) = p(\eta)\prod_{j=0}^{k-1}G_{A_j}(\eta_j, \eta_j).\]
\end{prop}

Suppose $V = \{x_1, x_2, \ldots, x_k\} \subset A$ and let $A_j = A \backslash \{x_1, \ldots, x_{j-1}\}$. Let
\[F(A;  x_1, x_2, \ldots, x_k) = \prod_{j=1}^{k}G_{A_j}(x_j, x_j).\]
The next lemma shows that it does not depend on the order we write the vertices $x_1, x_2, \ldots, x_k$.
\begin{lem} \label{lem3.8}
    If $\sigma : \{1, ..., k\} \rightarrow \{1, ..., k\}$ is a permutation, then
    \[F(A;x_{\sigma(1)}, ..., x_{\sigma(k)}) = F(A; x_1,...,x_k).\]
\end{lem}

\begin{proof}
    This is trivial if $k=1$. The proof for $k = 2$ is omitted and can be found in the proof of Lemma 2.5 in \cite{ams_RE}. For general $k$, we want to establish the result for any permutation $\sigma$ that just interchanges two adjacent indices:
    \begin{equation*}
        \begin{split}
            & F(A; x_1, \ldots, x_k)\\
            & = F(A; x_1, \ldots, x_i, x_{i+1}, \ldots, x_k)\\
            & = G_A(x_1, x_1) \cdots G_{A_i}(x_i, x_i)G_{A_{i+1}}(x_{i+1}, x_{i+1})\cdots G_{A_k}(x_k, x_k)\\
            & = G_A(x_1, x_1) \cdots F(A_i; x_i, x_{i+1}) \cdots G_{A_k}(x_k, x_k)\\
            & = G_A(x_1, x_1) \cdots F(A_i; x_{i+1}, x_{i}) \cdots G_{A_k}(x_k, x_k)\\
            & = G_A(x_1, x_1) \cdots G_{A\backslash\{x_1,\ldots,x_{i-1}\}}(x_{i+1}, x_{i+1})G_{A\backslash\{x_1,\ldots,x_{i-1},x_{i+1}\}}(x_{i}, x_{i}) \cdots G_{A_k}(x_k, x_k)\\
            & = F(A; x_1, \ldots, x_{i+1}, x_i, \ldots, x_k).
        \end{split}
    \end{equation*}
    Since any permutation $\{x_{\sigma(1)}, \ldots, x_{\sigma(k)}\}$ can be achieved by ``swapping'' two adjacent indices of $\{x_1,\ldots,x_k\}$, $F(A;x_{\sigma(1)}, \ldots, x_{\sigma(k)}) = F(A; x_1,\ldots,x_k)$.
\end{proof}

\begin{defn}\label{defn3.9}
    If $V = \{x_1, x_2, \ldots, x_k\} \subset A$, let
    \begin{equation}
        F_V(A) = \prod_{j=1}^{k}G_{A_j}(x_j, x_j)
    \end{equation}
    where $A_j = A \backslash \{x_1, \ldots, x_{j-1}\}$. Recall \autoref{lem3.8}. This quantity is independent of the
    ordering of the vertices. We also make the following conventions.
    \begin{itemize}
        \item If $V = A$, we write just $F(A)$ for $F_A(A)$.
        \item If $V \not\subset A$, then $F_V(A) = F_{V\cap A}(A)$.
        \item If $\eta$ is a path, $F_\eta(A)=F_V(A)$ where $V$ is the set of vertices visited by $\eta$.
    \end{itemize}
\end{defn}

The next proposition restates \autoref{prop3.7}.
\begin{prop}\label{prop3.11}
    If $\eta = [\eta_0, \eta_1, \ldots, \eta_k] \in \mathcal{R}_A(x, \partial A)$, then
    \[\hat{p}_A(\eta) = p(\eta)F_\eta(A).\]
\end{prop}

\subsection{Determinant of the Laplacian}

\begin{prop}\label{prop3.12}
    If $V \subset A$, then
    \[F_V(A) = \frac{\det G_A}{\det G_{A\backslash V}} = \frac{\det \mathcal{L}_{A\backslash V}}{\det \mathcal{L}}= \det \tilde{G}_V,\]
    where $\tilde{G}_V$ is the matrix $G_A$ restricted to the rows and columns indexed by $V$. In particular,
    \[F(A) = \frac{1}{\det \mathcal{L}}= \det G_A.\]
    Here $\mathcal{L} = \mathcal{L}_A = G^{-1}_A$.
\end{prop}

\begin{rem}
    Recall that $G_V = (I-P_V)^{-1}$ where $P_V$ is the transition matrix restricted to the rows and columns indexed by $V$. This is not the same at $\tilde{G}_V$ which is defined as the matrix $G_A$ restricted to those rows and columns.
\end{rem}

\begin{prop}
    Suppose $X_j,j = 0, 1,\ldots$, is an irreducible Markov chain on $\overline{A} = A \cup \partial A$ and suppose that either the chain is transient or $\partial A \neq \emptyset$. Let $x, y$ be distinct points in $A$, and let $q(x, y)$ be the probability that the chain starting at $x$ reaches $y$ before leaving $A$ or returning to $x$. Then
    \[G_A(x, y) = q(x, y)F_V(A),\]
    where $V = \{x, y\}$.
\end{prop}

\begin{proof}
    Recall that
    \[G_A(x, y) = \sum_{n=0}^{\infty}p_n(x, y) = \sum_{\omega: x \rightarrow y} p(\omega).\]
    We can decompose $\omega$ into three parts: loop(s) rooted at $x$ in $A$, a path starting at $x$ that reaches $y$ before leaving A or returning to $x$, and loop(s) rooted at $y$ in $A\backslash\{x\}$.  Then $\omega$ can be written as
    \[\omega = G_A(x, x)q(x, y)G_{A\backslash\{x\}}(y, y)\]
    and hence
    \[\omega = q(x,y)F(A; x, y) = F_V(A).\]
\end{proof}

\begin{prop}
    Suppose $X_j,j = 0, 1,\ldots$, is an irreducible Markov chain on $\overline{A} = A \cup \partial A$ where $A = \{x_1, ..., x_n\}$ is finite and $\partial A \neq \emptyset$. Let $P_A = [p(x, y)]_{x,y \in A}$ denote the transition matrix restricted to $A$ and $G_A = (I - P_A)^{-1}$ the Green's function. Suppose $V = \{x_1, ..., x_k\}$. We will consider the Markov chain that corresponds to the original chain “viewed only when visiting points in $V$”. The transition matrix $\tilde{P}  = [\tilde{p}(x, y)]_{x,y \in V}$ is given by
    \[\tilde{p} = \sum_{\omega: x \rightarrow y}p(\omega),\]
    where the sum is over all paths $\omega = [\omega_0,...,\omega_r]$ with $r \geq 1$; $\omega_0 = x,\omega_r =y$; and $\omega_1,...,\omega_{r-1} \in A\backslash V$. Then\\
    (1)the Green's function $\tilde{G}_V := (I-\tilde{P})^{-1}$ is the same as $G_A$ restricted to rows and columns indexed by vertices in $V$ and\\
    (2)\[F_V(A) = \det \tilde{G}_V.\]
\end{prop}

\begin{proof}
    (1) Recall that $G_A = I + P_A + P^2_A + P^3_A + \cdots$. Then $\tilde{G}_V = I + \tilde{P}_A + \tilde{P}^2_A + \tilde{P}^3_A + \cdots$. $\tilde{P}^i_A$ denotes the probability that the $i^{th}$ return to $V$ is at $y$. Let $T_i = \min \{k > T_{i-1}\mid X_k \in V\}$, where $T_0 = 0$ and $i = 0, 1, 2, 3, \ldots $ Then
    \[\tilde{P}^i_A(x, y) = \mathbb{P}^x\{X_{T_i}=y\}\]
    and hence
    \begin{equation*}
        \begin{split}
            \tilde{G}_V & = \sum_{i=0}^{\infty} \mathbb{P}^x\{X_{T_i}=y\}\\
            & = \sum_{i=0}^{\infty} \sum_{k=0}^{\infty} \mathbb{P}^x\{X_k=y, T_i = k\}\\
            & = \sum_{k=0}^{\infty} \sum_{i=0}^{\infty} \mathbb{P}^x\{T_i = k, X_k=y\}\\
            & = \sum_{k=0}^{\infty} \mathbb{P}^x\{X_k=y\}\\
            & = G_A(x, y)
        \end{split}
    \end{equation*}
    where $x, y \in V$.

    (2) From (1), observe that $G_A(y, y) = \tilde{G}_V(y, y)$ and $G_{A\backslash\{y\}}(z, z) = \tilde{G}_{V\backslash\{y\}}(z, z)$ for all $y, z \in V$. Then
    \begin{equation*}
            F_V(A) = \prod_{j=1}^{k}G_{A_j}(x_j, x_j) = \prod_{j=1}^{k}G_{\{x_j, \ldots, x_n\}}(x_j, x_j) = \prod_{j=1}^{k}\tilde{G}_{\{x_j, \ldots, x_k\}}(x_j, x_j) = \det \tilde{G}_V.
    \end{equation*}
\end{proof}

\subsection{Laplacian random walk}
\hfill\\

In this section, ``loop-erased walk from $x$ to $V \subset \partial A$'' is interpreted as follows.
\begin{itemize}
    \item Run the Markov chain starting at $x$ stopped at time $T$, the first time that $X_T \in \partial A$.
    \item Condition on the event that $X_T \in V$.
    \item Erase loops.
\end{itemize}
We can view this as a process $\hat{X}_0, \hat{X}_1, \ldots, \hat{X}_T$ where $\hat{X}_1, \ldots, \hat{X}_{T-1} \in A$ and $\hat{X}_T \in V \subset \partial A$. This chain does \textit{not} ``return'', so the conditional distribution of $\hat{X}_k$ given $\hat{X}_0, \ldots, \hat{X}_{k-1}$ depends on the entire past and not just on the value $\hat{X}_{k-1}$ and hence this is \textit{not} a Markov process. 

For $z \in \overline{A}$ and $w \in \partial A$, $\mathcal{K}_A(z, w)$ is the set of paths $\omega = [\omega_0, \omega_1, \ldots, \omega_k]$ for some $k$ with $\omega_0=z, \omega_k=w$ and $\omega_1, \ldots, \omega_{k-1} \in A$. If $x \in A$, we will write
\[\mathcal{K}_{A_x}(x, w)\]
where $A_x = A \backslash \{x\}$. That is, we turn $x$ into a boundary point.
\begin{defn}
    Suppose $z,w$ are distinct points in $\partial A$. Then the \textit{boundary Poisson kernel}, denoted $H_{\partial A}(z, w)$, is the measure of the set of paths starting at $z$, ending at $w$, and otherwise staying in $A$. If $z=w$ we define $H_{\partial A}(z,z)=1$. If $V \subset \partial A$, we write
    \[H_{\partial A}(z, V) = \sum_{w \in V}H_{\partial A}(z, w).\]
\end{defn}

\begin{prop}
    If $x \in A$, $z \in \partial A$,
    \[H_A(x, z) = G_A(x, x)H_{\partial A_x}(x, z),\]
    \[H_{\partial A_x}(x, z) = \sum_{y \in \overline{A}}p(x, y)H_{A_x}(y, z),\]
    where $A_x = A \backslash \{x\}.$
    More generally, if $V \subset \partial A$,
    \[H_A(x, V) = G_A(x, x)H_{\partial A_x}(x, V).\]
\end{prop}

\begin{proof}
    For $x \in A$ and $z \in \partial A$, note that
    \[H_A(x, z) = \sum_{\omega \in \mathcal{K}_{A}(x, z)}p(\omega).\]
    We can write $\omega$ as
    \[\omega = l_0 \oplus \omega^\ast\]
    with $l_0 \in \mathcal{K}_{A}(x, x)$, $\omega^\ast \in \mathcal{K}_{A_x}(x, z)$. Then
    \begin{equation*}
            H_A(x, z) = H_{\partial A_x}(x, z)\sum_{n=0}^{\infty}p_n(x, x) = H_{\partial A_x}(x, z)G_A(x, x).
    \end{equation*}
    To prove the second equality, let $T_{\partial A_x} = \min\{k \geq 1: X_k \in \partial A_x\}$, which denotes the first index that the chain visits $\partial A_x$. Then
    \begin{equation*}
        \begin{split}
            H_{\partial A_x}(x, z) & = \mathbb{P}^x\{X_{T_{\partial A_x}}=z\}\\
            & = \sum_{y\in \overline{A}} \mathbb{P}^x\{X_1 = y, X_{T_{\partial A_x}}=z\}\\
            & = \sum_{y\in \overline{A}} \mathbb{P}^x\{X_1 = y\} \mathbb{P}^x\{X_{T_{\partial A_x}}=z \mid X_1 = y\}\\
            & = \sum_{y\in \overline{A}} p(x, y)H_{A_x}(y, z).
        \end{split}
    \end{equation*}
    More generally, if $V \subset \partial A$,
    \begin{equation*}
        H_A(x, V) = \sum_{z \in V} H_A(x, z) = G_A(x, x)\sum_{z \in V}H_{\partial A_x}(x, z) = G_A(x, x)H_{\partial A_x}(x, V).
    \end{equation*}
\end{proof}

\begin{lem}
    Suppose $x \in A$, $V \subset \partial A$ and $H_A(x,V)>0$. Let $\{X_n\}$ be the Markov chain started at $x$ and $T =\min\{k\geq0:X_k \in \partial A\}$. Then for $y \in A_x \cup V$,
    \begin{equation*}
        \mathbb{P}^x\{X_1=y \mid X_T \in V, x \notin \{X_1, \ldots, X_{T-1}\}\} = \frac{p(x,y)H_{A_x}(y, V)}{H_{\partial A_x}(x, V)}.
    \end{equation*}
    Here $A_x = A \backslash \{x\}$.
\end{lem}

\begin{proof}
    We prove this lemma using conditional probability:
    \begin{equation*}
        \begin{split}
            & \mathbb{P}^x\{X_1=y \mid X_T \in V, x \notin \{X_1, \ldots, X_{T-1}\}\}\\
            & = \frac{\mathbb{P}^x\{X_1=y, X_T \in V, x \notin \{X_1, \ldots, X_{T-1}\}\}}{\mathbb{P}^x\{X_T \in V, x \notin \{X_1, \ldots, X_{T-1}\}\}}\\
            & = \frac{\mathbb{P}^x\{X_1=y, X_T \in V, x \notin \{X_1, \ldots, X_{T-1}\}\}}{H_{\partial A_x}(x, V)}\\
            & = \frac{\mathbb{P}^x\{X_1=y\}\mathbb{P}^x\{X_T \in V, x \notin \{X_1, \ldots, X_{T-1}\} \mid X_1=y\}}{H_{\partial A_x}(x, V)}\\
            & = \frac{p(x,y)H_{A_x}(y, V)}{H_{\partial A_x}(x, V)}.
        \end{split}
    \end{equation*}
\end{proof}

\subsection{Putting the loops back on the path}
\hfill\\

In this section, we consider the joint distribution of the LERW and the loops erased. A path $\omega \in \mathcal{K}_A(x, \partial A)$ is decomposed into its loop erasure $\eta$ and a collection of loops $\boldsymbol{l}=\{\ell_0, \ell_1, \ldots, \ell_{k-1}\}$ where $k = |\eta|$ is the length of $\eta$. The loop $\ell_j$ is an element of $\mathcal{K}_{A_j}(\eta_j, \eta_j)$. We write the joint distribution as
\[p(\eta, \boldsymbol{l}) = p(\eta)p(\boldsymbol{l}) \quad \text{where} \quad p(\boldsymbol{l}) = \prod_{j=0}^{k-1}p(\ell_j).\]
Recall \autoref{prop3.11}. We can write
\[p(\eta, \boldsymbol{l}) = \hat{p}(\eta)p(\boldsymbol{l})F_{\eta}(A)^{-1}.\]

\begin{prop}
    Consider LERW from $x$ to $\partial A$, that is, we assign each SAW $\eta = [\eta_0,\ldots,\eta_k]$ from $x$ to $\partial A$ in $A$ probability $\hat{p}(\eta)$. Suppose that we choose conditionally independent loops $\boldsymbol{l}$ given $\eta$ with probability
    \[F_\eta(A)^{-1}p(\boldsymbol{l}) = F_\eta(A)^{-1}p(\ell_0)\ldots p(\ell_{k-1}),\]
    assuming $\ell_j \in \mathcal{K}_{A_j}(\eta_{j-1}, \eta_{j-1})$. Then if we put the loops on the curve and concatenate as in the proof of \autoref{prop3.7}, the path we get has the distribution of the Markov chain starting at $x$ stopped at $\partial A$.

\end{prop}

\begin{prop}
    The following gives a way to sample the loops $\ell_0, \ldots, \ell_{k-1}$ given $\eta$.
    \begin{itemize}
        \item Start independent Markov chains $X^j$ at each $\eta_j$ and let \[T^j = \min \{n \geq 0: X^j_n \in \partial A_j\}.\]
        \item Let \[\rho^j = \max \{m<T^j: X^j_m = x_j\}.\]
        \item Output \[\ell_j = [X^j_0, X^j_1, \ldots, X^j_{\rho^j}].\]
    \end{itemize}
\end{prop}

\begin{proof}
    We want to show the probability that we choose conditionally independent loops $\boldsymbol{l}$ given $\eta$ is
    \[F_\eta(A)^{-1}p(\boldsymbol{l}) = F_\eta(A)^{-1}p(\ell_0)\ldots p(\ell_{k-1}),\]
    where $\ell_j \in \mathcal{K}_{A_j}(\eta_{j-1}, \eta_{j-1})$. Note that
    \[\mathbb{P} \{\ell_0, \ldots, \ell_{k-1} \mid \eta = [\eta_0, \ldots \eta_k]\} = \prod_{j=0}^{k-1} \mathbb{P} \{\ell_j \mid \eta\}\]
    and
    \begin{equation*}
        \begin{split}
            \mathbb{P} \{\ell_j \mid \eta\} & = \sum_{\omega_j \in \mathcal{K}_A(\eta_j, \partial A_j)}p(\omega_j)\\
            & = \sum_{\hat{\omega_j} \in \mathcal{K}_{A\backslash \{\eta^j\}}(\eta_j, \partial A_j)} p(\ell_j \oplus \hat{\omega_j})\\
            & = p(\ell_j)\sum_{\hat{\omega_j} \in \mathcal{K}_{A\backslash \{\eta_j\}}(\eta_j, \partial A_j)} p(\hat{\omega_j})\\
            & = p(\ell_j)H_{\partial A_{j+1}}(\eta_j, \partial A_j).
        \end{split}
    \end{equation*}
    Then we substitute back into $\mathbb{P} \{\ell_0, \ldots, \ell_{k-1} \mid \eta\}$:
    \begin{equation*}
        \begin{split}
            \mathbb{P} \{\ell_0, \ldots, \ell_{k-1} \mid \eta\} & = \prod_{j=0}^{k-1} p(\ell_j)H_{\partial A_{j+1}}(\eta_j, \partial A_j)\\
            & = \prod_{j=0}^{k-1} \frac{1}{G_{A_j}(\eta_j, \eta_j)} p(\ell_0)p(\ell_1) \cdots p(\ell_{k-1})\\
            & = F_\eta(A)^{-1}p(\boldsymbol{l}).
        \end{split}
    \end{equation*}
\end{proof}

\subsection{Wilson's algorithm}
\hfill\\

\begin{defn}
    \hfill\\
    \begin{itemize}
        \item A \textit{tree} is a simple undirected graph such that for any two distinct vertices $x, y$ there is exactly one SAW in the graph with initial vertex $x$ and terminal vertex $y$.
        \item A \textit{spanning tree} of a graph is a subgraph containing all the vertices that is a tree.
    \end{itemize}
\end{defn}

Any tree on $n$ vertices must have exactly $n-1$ edges. A finite graph only has a finite number of spanning trees and hence it makes sense to choose one ``at random''.

\begin{defn}
    A \textit{uniform spanning tree} (UST) of a graph is a random spanning tree chosen from the uniform distribution on all spanning trees.
\end{defn}

The word “uniform” in uniform spanning tree refers to the probability distribution on trees and shows that it uniform over all trees. The algorithm found by David Wilson \cite{Wilson} for choosing a spanning tree selects from the uniform distribution, and it also allows us to get an expression for the number of trees.

Assume a graph has $n+1$ vertices. Choose any ordering of the vertices and write them as $V = \{x_0, x_1, \ldots, x_n\}$. Let $V_0 = \{x_0\}$. Start simple random walk at $x_1$ and stop it when it reaches $x_0$. Erase loops chronologically from the path. Let $E_1$ be the set of edges from the loop erasure of this path and let $V_1$ be the set of vertices that are adjacent to an edge in $E_1$. Recursively, if $V_k = V$ stop. Otherwise, let $x_j$ be the vertex of smallest index not in $V_k$. Take random walk starting at $x_j$ stopped when it reaches $V_k$. Erase loops, and add the remaining edges to $E_k$ giving $E_{k+1}$ and let $V_{k+1}$ be the set of vertices that are adjacent to an edge in $E_{k+1}$. We call this \textit{Wilson's algorithm} of generating random spanning trees.

\begin{prop}
    If $\mathcal{T}$ is a spanning tree of a connected graph with vertices $V = \{x_0, \ldots, x_n\}$, then the probability that $\mathcal{T}$ is chosen in Wilson's algorithm is
    \[\left[ \prod_{j=1}^{n} \deg (x_j)\right]^{-1} F(A)\]
    where $A = \{x_1, \ldots, x_n\}.$ Recall \autoref{prop3.12}, which states that $F(A) = \det G_A$. Note that the probability does not depend on how the points were ordered or on the tree $\mathcal{T}$.
\end{prop}

Since we are choosing from a uniform distribution and we know the probability of picking a particular element, we know how many elements are in the set which gives the following corollary.

\begin{cor}\label{Cor3.24}
    The total number of spanning trees is
    \[\left[ \prod_{j=1}^{n} \deg (x_j)\right] F(A)^{-1} = \left[ \prod_{j=1}^{n} \deg (x_j)\right] \det \mathcal{L}_A.\]
\end{cor}

\begin{prop}
    The complete graph on n+1 vertices $\{x_0, x_1, \ldots, x_n\}$ is the (undirected) simple graph with the maximal number of edges. That is, every pair of distinct vertices is connected by an edge. In particular, each vertex has degree $n$. The number of spanning trees of the complete graph on $n + 1$ vertices is
    \[(n+1)^{n-1}.\]    
\end{prop}

\begin{proof}
    We first compute
    \[F(A) = \prod_{j=1}^{n} G_{A_j}(x_j, x_j), \quad A_j = \{x_j, x_{j+1}, \ldots, x_n\}.\]
    For vertices $\{x_0, x_1, \ldots, x_n\}$, suppose the chain starts at $x_j$. The chain is ``killed'' if the chain leaves $A_j$. Observe that the probability that $x_j$ leaves $A_j$ is
    \[\frac{j}{n}+\frac{n-j}{n}\frac{j}{j+1}.\]
    This is equivalent to the probability that the chain never returns to $x_j$ starting at $x_j$.
    Recall that $G_{A_j}(x_j, x_j)$ is the expected value of a geometric distribution:
    \[G_{A_j}(x_j, x_j) = \frac{1}{\frac{j}{n}+\frac{n-j}{n}\frac{j}{j+1}} = \frac{n(j+1)}{j(n+1)}.\]
    Then we have
    \[F(A) = \prod_{j=1}^{n} \frac{n(j+1)}{j(n+1)} = \frac{n^n}{(n+1)^{n-1}},\]
    and using \autoref{Cor3.24}, the total number of spanning trees is
    \[n^n\frac{(n+1)^{n-1}}{n^n} = (n+1)^{n-1}.\]
\end{proof}

\section{Loop Soups}

We focus on the loops in this chapter. We define a growing loop at a point and the growing loop configuration in $A$. Then we show how they can be constructed in three different ways: (ordered) growing loops, rooted loops, and unrooted loops.

\subsection{Growing loop at a point}
\hfill\\

In this section, we give a model for a ``growing loop'' rooted at a vertex $x$ in $A$. We let $l_t$ denote a continuous time Markov chain for time $t \geq 0$ with state space $\mathcal{K}: = \mathcal{K}(x) : = \mathcal{K}_A(x, x)$ starting at the trivial path. We allow the trivial loop. We assume that $l_0$ is the trivial loop, and the process grows at the end. That is to say, if $s \leq t$, $l_s$ is an initial segment of $l_t$, that is $l_t = l_s \oplus l^*$ for another loop $l^*$. The distribution at time $t = 1$ will be that of the loops erased in the chronological loop-erasing procedure.

For $n > 0$, let $\mathcal{K}_n = \mathcal{K}_n(x)$ be the set of loops starting at $x$ that return to $x$ exactly $n$ times. We call the loops in $\mathcal{K}_1$ elementary loops and note that each loop in $\mathcal{K}_k$ can be written uniquely as
\begin{equation}\label{(4.1)}
    l = l^1 \oplus l^2 \oplus \cdots \oplus l^k, \quad l^j \in \mathcal{K}_1.
\end{equation}
The total measure of the set of elementary loops at $x$ is denoted by
\[f_x: = p(\mathcal{K}_1) = \sum_{l \in \mathcal{K}_1}p(l) = 1 - \frac{1}{G_A(x, x)} = \frac{G_A(x, x)-1}{G_A(x, x)}.\]
If $p$ is the transition matrix of a Markov chain, $f_x$ is the probability that the chain starting at $x$ returns to $x$ without leaving $A$. Let $\nu$ denote the probability measure on $\mathcal{K}_1$ given by $\nu(l) = p(l)/f_x$. Using \autoref{(4.1)} we see that $p[\mathcal{K}_k]=f^k_x$.

\begin{defn}
    (Growing loop)
    \begin{itemize}
        \item Let $K_t$ be the number of elementary loops in the growing loop at time $t$. We assume that $K_t$ is a continuous time Markov chain with state space $\mathbb{N}$ and rate \[\lambda(n, n+k) = \frac{1}{k}f^k_x.\]
        \item Choose $l^1, l^2, \ldots$ to be independent, identically distributed loops in $\mathcal{K}_1$ with distribution $\nu$.
        \item Let \[l_t = l^1 \oplus l^2 \oplus \cdots \oplus l^{K_t},\]where the right-hand side is defined to be the trivial loop if
        $K_t=0$.
    \end{itemize}
\end{defn}

\begin{defn}
    (An alternative definition for growing loop)
    \begin{itemize}
        \item $l_t$ is a continuous time Markov chain taking values in $\mathcal{K}$ with $l_0$ being the trivial loop and rates \[\lambda(l, l \oplus l') = \frac{1}{k}p(l'), \quad l' \in \mathcal{K}_k.\]
    \end{itemize}
\end{defn}

The process $K_t$ is a \textit{negative binomial process} with parameter $G_A(x, x)^{-1}=1-f_x$. $K_t$ is counting $k$ elementary loops (failures), given $t$ chains starting at $x$ either never return to $x$ or leave $A$ before returning to $x$ (successes). From this we see that
\[\mathbb{P}\{K_t=k\}=\frac{1}{G_A(x, x)^t}\frac{\Gamma(k+t)}{k!\Gamma(t)}f^k_x, \quad k = 0, 1, 2, \ldots,\]
and hence
\begin{equation}\label{(4.4)}
    \mathbb{P}\{l_t=l\}=\frac{1}{G_A(x, x)^t}\frac{\Gamma(k+t)}{k!\Gamma(t)}p(l) \quad \text{if} \quad l \in \mathcal{K}_k.
\end{equation}
When $t=1$, we get the simpler expression
\begin{equation}\label{(4.5)}
    \mathbb{P}\{l_1=l\}=\frac{1}{G_A(x, x)^t}p(l).
\end{equation}

\begin{exmp} \label{exmp 4.1}
    Suppose we have a Markov chain on the state space $\overline{A} = \{1, 2, 3\}$ with transition probabilities
    \[p(1, 3) = p(2, 3) = \frac{1}{3}, \quad p(3, 3)=1\]
    \[p(1, 1) = \frac{1}{3}, \quad p(1, 2) = \frac{1}{3}, \quad p(2, 1) = \frac{1}{2}, \quad p(2, 2) = \frac{1}{6}.\]
    Let $A = \{1, 2\}$. We will find a number of quantities with $x = 1$ and with $x=2$.

    
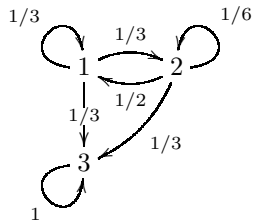
\[\xymatrix{
        1 \ar[d]|{1/3} \ar@/^/[r]^-{1/3} \ar@(l,u)[]^-{1/3} & 2 \ar@/^/[l]^-{1/2} \ar@/^/[dl]^-{1/3} \ar@(r,u)[]_{1/6}\\
        3 \ar@(l,d)[]_{1}}\]

    (1) $f_x$. Observe that the set of elementary loops at $x$ is of the form \[\{[x, x], [x, y, x], [x, y, y, x], [x, y, y, y, x], \ldots\}.\] Then for $x=1$
    \[f_1 = \frac{1}{3}+\frac{1}{3}\left(\frac{1}{2}\right) + \frac{1}{3}\left(\frac{1}{6}\right)\left(\frac{1}{2}\right) + \frac{1}{3}\left(\frac{1}{6}\right)^2\left(\frac{1}{2}\right) + \ldots = \frac{8}{15}\]
    and for $x = 2$
    \[f_2 = \frac{1}{6}+\frac{1}{2}\left(\frac{1}{3}\right) + \frac{1}{2}\left(\frac{1}{3}\right)\left(\frac{1}{3}\right) + \frac{1}{2}\left(\frac{1}{3}\right)^2\left(\frac{1}{3}\right) + \ldots = \frac{5}{12}.\]
    
    (2) $G_A(x, x)$. Recall that $f_x = 1 - \frac{1}{G_A(x, x)}$. Then we have
    \[G_A(1, 1) = \frac{1}{1-f_1} = \frac{15}{7}, \quad G_A(2, 2) = \frac{1}{1-f_2} = \frac{12}{7}.\]

    (3) $\mathbb{P}\{l_1 = [x]\}.$ Using \autoref{(4.5)}, the probability that the growing loop at time $t = 1$ is trivial is
    \[\mathbb{P}\{l_1 = [x]\} = \frac{1}{G_A(x, x)}p([x]) = \frac{1}{G_A(x, x)}\]
    and hence
    \[\mathbb{P}\{l_1 = [1]\} = \frac{7}{15}, \quad \mathbb{P}\{l_1 = [2]\}  = \frac{7}{12}.\]

    (4) $\mathbb{P}\{l_1 \text{ is of length 2}\}.$ Observe that $\mathbb{P}\{l_1 \text{ is of length 2}\} = \mathbb{P}\{l_1 = [x, x, x]\} + \mathbb{P}\{l_1 = [x, y, x]\}$. Use \autoref{(4.4)}. Then for $x=1$
    \[\mathbb{P}\{l_1 \text{ is of length 2}\} = \frac{7}{54}\]
    and for $x=2$
    \[\mathbb{P}\{l_1 \text{ is of length 2}\}= \frac{49}{432}.\]

    \begin{multline*}
        (5)\text{ }
        \mathbb{P}\{l_1 = [x, x, x] \mid l_1 \text{ is of length 2}\} = \frac{\mathbb{P}\{l_1 = [x, x, x]\}}{\mathbb{P}\{l_1 \text{ is of length 2}\}}.\\
    \end{multline*}
    We have calculated all the quantities.

    (6) Let $\rho(t)$ be the probability that the growing loop at time $t$ is of length 2 and let $\tilde{\rho}(t)$ be the probability that the growing loop at time $t$ is $[x, x, x]$. Then we have
    \[\lim_{t \downarrow 0}\frac{\tilde{\rho}(t)}{\rho(t)} = \frac{1}{4}, \quad \lim_{t \rightarrow \infty}\frac{\tilde{\rho}(t)}{\rho(t)} = 1,\]
    which can be seen through simplification.
\end{exmp}

\subsection{Growing loop configuration in $A$}
\hfill\\

The growing loop configuration in $A$ is a collection of growing loops at each vertex. It depends on an ordering $\sigma$ of the vertices, and, as before, we let $A_j = A \backslash \{x_1,\ldots,x_{j-1}\}$. Let us write $\mathcal{G} = \mathcal{G}_{A, \sigma}$ to denote the set of $n$-tuples of loops
\[\boldsymbol{l} = (l^1, \ldots, l^n)\]
such that $l^j$ is a loop in $A_j$ rooted at $x_j$.

\begin{defn}
    Given the ordering of $A$, the growing loop configuration is an $n$-tuple
    \[\boldsymbol{l}_t = (l^1_t, l^2_t, \ldots, l^n_t)\]
    where $l^1_t, l^2_t, \ldots, l^n_t$ are independent and $l^j_t$ is a growing loop at $x_j$ in $A_j$ as in the previous section. 
\end{defn}

\begin{prop}\label{prop4.8}
    If $\sigma$ is an ordering of $A$ and $p$ is an integrable nonnegative weight, then for $\boldsymbol{l} \in \mathcal{G}$,
    \[\mathbf{P}\{\boldsymbol{l}_t = \boldsymbol{l}\} = \frac{p(\boldsymbol{l})}{(\det G)^t}\prod_{i=1}^{n}\frac{\Gamma(j_i+t)}{j_i!\Gamma(t)}\]
    where $j_i$ is the number of times that loop $l_i$ returns to the origin as in the previous section and $p(\boldsymbol{l}) = p(l^1)p(l^2)\cdots p(l^n)$.

\end{prop}

We view the growing loop configuration as a \textit{``loop soup''}. The growing loop configuration depends on the ordering of vertices and we  write $\sigma$ to denote the particular ordering of the vertices. Let \[\mathcal{J}_\sigma = \bigcup_{j=1}^{n}\mathcal{K}_{A_j}(x_j, x_j).\]
For each $l \in \mathcal{K}_{A_j}(x_j, x_j)$, let $\beta(l)$ denote the number of times that the loop returns to $x_j$. For each $l \in \mathcal{J}_\sigma$, let $N_t(l)$ denote the number of times that $l$ has been added onto a loop by time $t$. Then $N_t(l): = \{N_t(l): l \in \mathcal{J}_\sigma\}$ are independent Poisson processes with parameters
\[\lambda_l = \frac{1}{\beta(l)}p(l), \quad \text{if } l \in \mathcal{K}_{A_j}(x_j, x_j).\]

\begin{exmp}
    Consider the Markov chain in \autoref{exmp 4.1} and suppose that the vertices are ordered $1, 2$ and $l_t = (l^1_t, l^2_t)$ is the corresponding growing loop configuration.

    (1) For each $t$, we will find the probability that the loop $l^1_t$ visits the vertex 2 using two different approaches. Loops that visit the vertex 2 can be separated based on the number of times that the loop returns to vertex 1. Loops that visit the vertex 2 $k$ times are Poisson processes with parameter
    \[\frac{1}{k}\left[f_1^k - \left(\frac{1}{3}\right)^k\right].\]
    Since $N_t$ are independent Poisson processes, the number of times that loops that visit the vertex 2 has been added by time $t$ has parameter
    \begin{equation*}
            \sum_{k=1}^{\infty} \frac{1}{k}\left[f_1^k - \left(\frac{1}{3}\right)^k\right] = \sum_{k=1}^{\infty} \frac{1}{k}\left[\left(\frac{8}{15}\right)^k - \left(\frac{1}{3}\right)^k\right] = \log \frac{10}{7}.
    \end{equation*}
    Then we have
    \[\mathbb{P}\{l^1_t \text{ visits the vertex 2}\} = 1 - e^{-(\log \frac{10}{7})t} = 1 - \left(\frac{7}{10}\right)^t.\]
    The other approach is
    \begin{equation*}
        \begin{split}
            \mathbb{P}\{l^1_t \text{ visits the vertex 2}\} & = 1 - \mathbb{P}\{l^1_t \text{ doesn't visit the vertex 2}\}\\
            & = 1 - \frac{1}{G_A(1, 1)}\sum_{k=0}^{\infty}\frac{\Gamma(k+t)}{k!\Gamma(t)}\left(\frac{1}{3}\right)^k\\
            & = 1 - \left(\frac{7}{10}\right)^t.
        \end{split}
    \end{equation*}

    (2) For each t, we will find the expected number of times vertex 2 appears in $l^1_t$. Note that
    \[\mathbb{E}[\# \text{ 2's in } l^1_t] = \sum_{k=1}^{\infty}\mathbb{E}[\# \text{ 2's in } l^1_t \mid K_t = k] \cdot \mathbb{P}\{K_t = k\}.\]
    Then 
    \begin{equation*}
        \begin{split}
            \sum_{k=1}^{\infty}\mathbb{E}[\# \text{ 2's in } l^1_t \mid K_t = k] & = \sum_{l^j \in \mathcal{K}_1}\mathbb{P}\{l^1_t = l^1 \oplus l^2 \oplus \cdots \oplus l^k\}(\# \text{ 2's in } l^1_t)\\
            & = k \cdot \sum_{l \in \mathcal{K}_1} \frac{p(l)}{f_x}(\# \text{ 2's in } l)\\
            & = k \cdot \sum_{n=1}^{\infty} \frac{(1/6)^n}{8/15}n\\
            & = \frac{9}{20}k.
        \end{split}
    \end{equation*}
    The second equality holds since $l^j$ are independent. Recall that the process $K_t$ is a \textit{negative binomial process} and hence
    \begin{equation*}
            \mathbb{E}[\# \text{ 2's in } l^1_t] = \sum_{k=1}^{\infty} \frac{9}{20}k \cdot \mathbb{P}\{K_t = k\} = \frac{9}{20} \mathbb{E}[K_t] = \frac{18}{35}t.
    \end{equation*}
\end{exmp}

\begin{prop}
    Order the vertices $A = \{x_1, \ldots, x_n\}$ and define the measure on ordered loops $m^*$ as follows. Let $A_j = A\backslash\{x_1,...,x_{j-1}\}$. If $l$ is a loop in $A_j$ rooted at $x_j$, then
    \[m^*(l) = \frac{p(l)}{\beta(l)}\]
    where $\beta(l) = \# \{k:1\leq k \leq |l|, l_k = x_j\}$, and $m^*(l)=0$ for all other loops. Then
    \[\sum_{l} m^*(l) = \log \det G_A.\]
\end{prop}

\begin{proof}
    We first calculate loops in $A_j$ rooted at $x_j$:
    \begin{equation*}
        \sum_{l \in K_{A_j}(x_j, x_j)} m^*(l) = \sum_{l \in K_{A_j}(x_j, x_j)} \frac{p(l)}{\beta(l)} = \sum_{k=1}^{\infty} \frac{1}{k}f^k_{x_j} = - \log \left(\frac{1}{G_{A_j}(x_j, x_j)}\right).
    \end{equation*}
    Then it follows that
    \begin{equation*}
        \begin{split}
            \sum_{l} m^*(l) & = \sum_{j=1}^{n}- \log \left(\frac{1}{G_{A_j}(x_j, x_j)}\right)\\
            & = - \log \left(\frac{1}{\prod_{j=1}^{n}G_{A_j}(x_j, x_j)}\right)\\
            & = - \log \left(\frac{1}{F(A)}\right)\\
            & = - \log \left(\frac{1}{\det G_A}\right)\\
            & = \log \det G_A.
        \end{split}
    \end{equation*}
    The third equality uses \autoref{defn3.9} and the penultimate equality uses \autoref{prop3.12}.
\end{proof}

We now focus on $t=1$, for which \autoref{(4.5)} shows that the conclusion of \autoref{prop4.8} can be written as
\begin{equation}\label{(4.11)}
    \mathbf{P}\{\boldsymbol{l}_1 = \boldsymbol{l}\} = \frac{p(\boldsymbol{l})}{\det G}\prod_{i=1}^{n}\frac{p(\ell^i)}{G_{A_j}(x_j, x_j)}.
\end{equation}

\begin{prop}
    See Proposition 3.9 and Proposition 3.10 in \cite{ams_RE}. The propositions show that we can construct realizations of a Markov chain by starting with realizations of the loop-erased random walk or the uniform spanning tree, taking conditionally independent realizations of the loop soup at time $t = 1$ and combining them.
\end{prop}

\begin{proof}
    For Proposition 3.9 in \cite{ams_RE}, note that
    \[p(\omega) = p(\eta)\prod_{j=1}^{k-1}p(l^j).\]
    We want to show that
    \[p(\omega) = \hat{p}(\eta)\prod_{j=0}^{k-1}\mathbb{P}\{l^j_1 = l^j\}.\]
    Using \autoref{(4.11)} and \autoref{prop3.7}, we can see that
    \begin{equation*}
        \begin{split}
            \hat{p}(\eta)\prod_{j=0}^{k-1}\mathbb{P}\{l^j_1 = l^j\} & = \hat{p}(\eta)\prod_{j=0}^{k-1}\frac{p(l^j)}{G_{A_j}(\eta_j, \eta_j)}\\
            & = p(\eta) \prod_{j=0}^{k-1}G_{A_j}(\eta_j, \eta_j)\prod_{j=0}^{k-1}\frac{p(l^j)}{G_{A_j}(\eta_j, \eta_j)}\\
            & = p(\eta) \prod_{j=0}^{k-1}p(l^j).
        \end{split}
    \end{equation*}
    
    For Proposition 3.10 in \cite{ams_RE}, we want to show that
    \[\prod_{\omega} p(\omega) = \mathbb{P}\{\text{the tree } \mathcal{T} \text{ is chosen}\}\cdot \prod_{j=1}^{n} \mathbb{P}\{l^j_1 = l^j\}.\]
    Observe that the left-hand side is
    \[\prod_{\omega} p(\omega) = \prod_{e \in \mathcal{T}}p(e) \prod_{j=1}^{n}p(l^j) = p(\mathcal{T})\prod_{j=1}^{n}p(l^j).\]
    Using Proposition 2.27 in \cite{ams_RE} and \autoref{defn3.9}, we can see that the right-hand side is
    \begin{equation*}
        \begin{split}
            \mathbb{P}\{\text{the tree } \mathcal{T} \text{ is chosen}\} & \cdot \prod_{j=1}^{n} \mathbb{P}\{l^j_1 = l^j\} = p(\mathcal{T})F(A)\prod_{j=1}^{n} \mathbb{P}\{l^j_1 = l^j\}\\
            & = p(\mathcal{T}) \prod_{j=1}^{n}G_{A_j}(x_j, x_j) \prod_{j=1}^{n}\frac{p(l^j)}{G_{A_j}(\eta_j, \eta_j)}\\
            & = p(\mathcal{T}) \prod_{j=1}^{n}p(l^j).
        \end{split}
    \end{equation*}
\end{proof}

\subsection{Rooted loop soup}
\hfill\\

In the previous section, the ordered loop soup depends on an ordering $\sigma$ of the vertices. In this section, we explore a different measure and corresponding loop soup that does not depend on the ordering.

\begin{defn}
    Suppose $\omega = [\omega_0, \omega_1, \ldots, \omega_n]$ is a loop. We write $|\omega| = n$ for the number of steps in $\omega$. Then we define the translation $\tau$ as follows:
    \[\tau \omega = [\omega_1, \omega_2, \ldots, \omega_n, \omega_1].\]
\end{defn}
Note that $\tau \omega$ traverses the same loop as $\omega$ in the same direction but with a different starting point. We let $J(\omega)$ be the number of distinct loops obtained. Let $J(\omega;x)$ be the number of distinct loops obtained that are rooted at $x$. Then,
\begin{equation}
    J(\omega; x) = \frac{J(\omega)}{|\omega|}\# \{j: 1 \leq j \leq n, \quad \omega_j=x\}.
\end{equation}

\begin{defn}\label{defn4.15}
    The \textit{rooted loop measure} $\tilde{m}=\tilde{m}^p$ is the measure that assigns measure
    \begin{equation}
        \tilde{m}(l) = \frac{p(l)}{|l|}
    \end{equation}
    to every loop
    \[l \in \bigcup_{j=1}^{n}\mathcal{K}_A(x_j, x_j).\]
    The corresponding Poisson realization is called the \textit{rooted loop soup}.
\end{defn}
Note that the rooted loop measure is not a probability measure since
\begin{equation}\label{(4.17)}
    \sum_{l}\tilde{m}(l) = \log \det G_A.
\end{equation}

\begin{prop} \label{prop4.5}
    The following is a valid way to get the rooted loop soup.
    \begin{itemize}
        \item Choose an ordering of the vertices $\sigma$.
        \item Take a realization of the ordered loop soup with ordering $\sigma$.
        \item For each $\omega$ in the soup, choose a translation $\tau^k\omega$ where $k$ is chosen uniformly in $\{0, 1, \ldots, |\omega|-1\}$.
    \end{itemize}
\end{prop}

\begin{proof}
    For each $l \in \bigcup_{j=1}^{n}\mathcal{K}_A(x_j, x_j)$, we want to show that the \textit{rooted loop measure} of $l$ is
    \[\tilde{m}(l) = \frac{p(l)}{|l|}.\]
    Let $i = \min \{j: x_j \in l\}$, which denotes the minimum index in $l$. Recall that $J(l; x_i)$ denotes the number of distinct loops obtained by translation that are rooted at $x_i$. Let $\beta(l; x_i)$ denote the the number of times that $l$ returns to $x_i$. Since we are choosing a translation $\tau^k\omega$ uniformly, we get
    \begin{equation*}
            \tilde{m}(l) = \frac{1}{J(l)}\left[J(l;x_i)\frac{p(l)}{\beta(l;x_i)}\right] = \frac{1}{J(l)}\left[\frac{J(l)}{|l|}\beta(l;x_i)\frac{p(l)}{\beta(l;x_i)}\right] = \frac{p(l)}{|l|}.
    \end{equation*}
    $J(l)$ is the number of distinct loops obtained by translation and $\frac{1}{J(l)}$ shows that the translation is chosen uniformly.
\end{proof}

\subsection{(Unrooted) random walk loop measure}
\hfill\\

In the last section, we choose a loop and then randomize its starting point. We write $\ell$ for unrooted loops and $l$ or $\omega$ for rooted loops.

\begin{defn}\label{defn4.19}
    The \textit{unrooted loop measure} $m=m^p$ is the measure that assigns to each unrooted loop
    \begin{equation}
        m(\ell) = J(\ell)\frac{p(\ell)}{|\ell|}.
    \end{equation}
    The \textit{(unrooted) loop measure} is a Poissonian realization from $m$.
\end{defn}

Using \autoref{(4.17)}, we see that
\begin{equation}\label{(4.21)}
    \sum_{\ell}m(\ell) = \log \det G_A.
\end{equation}

\begin{prop}
    Each of the following methods is a valid way to get the unrooted loop soup.
    \begin{itemize}
        \item Take a realization directly from m.
        \item Take a realization of the rooted loop soup and then ``forget the root''.
        \item Choose an ordering $\sigma$, take a realization of the ordered loop soup with ordering $\sigma$, and then``forget the root''.
    \end{itemize}
\end{prop}

\begin{proof}
    (1) It directly follows from Definition 3.14.

    (2) Now, for each $\omega$ in the soup, distinct $\tau^k\omega$ are chosen together. From \autoref{prop4.5} we see that
    \[m(\ell) = J(\mathcal{\ell})\tilde{m}(\ell) = J(\mathcal{\ell})\frac{p(\ell)}{|\ell|},\]
    which corresponds to the \textit{unrooted loop measure} of $\ell$ in \autoref{defn4.19}.

    (3) We take a realization of the ordered loop soups, and then we repeat (2).
\end{proof}

\begin{prop}
    The probability that the unrooted loop soup has no loops in it at time $t$ is given by $(\det G_A)^{-t}$.
\end{prop}

\begin{proof}
    The \textit{(unrooted) loop soup} is a Poisson realization from $m$. Growing loops are independent Poisson processes with parameters
    \[m(\ell).\]
    Let $X$ be the number of loops in the unrooted loop soup at time $t$. Using \autoref{(4.21)}, we see that
    \[\lambda_{total} = \sum_{\ell}m(\ell) = \log \det G_A.\]
    Then we have $X \sim Poi(\log \det G_A)$ and hence
    \[\mathbb{P}\{X=0\} = \frac{(\lambda_{total}t)^0e^{-\lambda_{total}t}}{0!} = e^{-(\log \det G_A)t} = (\det G_A)^{-t}.\]
\end{proof}

\section*{Acknowledgments}  I would like to express my heartfelt gratitude to my mentor, Mark Olson, for his invaluable support during this research and for the time he took to guide me through the materials presented above. I am deeply indebted to Professor Gregory Lawler for holding Probability and Analysis lectures and first introducing me to LERW and Loop Soups. I would also like to thank Jinwoo Sung for his immense patience in answering my questions. A final thank you to Professor Peter May for organizing the REU; being extremely accessible, supportive and kind; and bringing together such an engaged community of students, mentors and speakers. Writing this paper requires a great deal of effort. Without any of these people, this paper would've either never been good or never been completed.

\end{document}